\documentclass[a4paper]{amsart}

\usepackage{amsmath,amsthm,amsfonts,amssymb,amscd}
\input{xy}
\xyoption{all}
\newtheorem{theo}{Theorem}[section]
\newtheorem{pro}[theo]{Proposition}
\newtheorem{lem}[theo]{Lemma}
\newtheorem{cor}[theo]{Corollary}

\newenvironment{sis}{\left\{\begin{aligned}}{\end{aligned}\right.}

\numberwithin{equation}{section}

\newcommand{\Z}{\mathbb{Z}}

\newcommand{\0}{\underline{0}}

\newcommand{\g}{\mathfrak{g}}
\newcommand{\h}{\mathfrak{h}}

\renewcommand{\d}{{\rm d}}

\newcommand{\Sq}{{\rm Sq}}

\newcommand{\ti}{\widetilde}

\begin{document}

\title{Deformations of the restricted Melikian Lie algebra}
\author{Filippo Viviani}

\date{07 March 2008}

\address{Humboldt Universit\"at zu Berlin, 
Institut f\"ur Mathematik, 
Rudower Chausse 25, 
10099 Berlin (Germany).}

\email{ viviani@math.hu-berlin.de}



\thanks{During the preparation of this paper, the author was supported by a grant from the 
Mittag-Leffler Institute of Stockholm.}

\begin{abstract}
We compute the infinitesimal deformations of the restricted Melikian Lie algebra
in characteristic $5$. 
\end{abstract}

\maketitle

\section{Introduction}

The restricted Melikian Lie algebra $M$ is a restricted simple Lie algebra of dimension 
$125$ defined over the prime field of characteristic $p=5$. It was introduced by 
Melikian \cite{MEL} and it is the only ``exceptional'' simple Lie algebra
which appears in the classification of restricted simple Lie algebras over a field
of characteristic $p\neq 2,3$ (see \cite{BW} for $p>7$ and \cite{PS3} for $p=5,7$).
The classification problem remains still open in characteristic $2$ and $3$, where 
several ``exceptional'' simple Lie algebras are known (see \cite[page 209]{STR}).

This paper is devoted to the study of the \emph{infinitesimal deformations}
of the restricted Melikian Lie algebra. The infinitesimal deformations have been 
computed for the other restricted simple Lie algebras in characteristic $p\geq 5$.
Rudakov (\cite{RUD}) proved that the simple Lie algebras of classical type are rigid,
in analogy of what happens in characteristic zero. On the other hand, 
the author computed the infinitesimal deformations of the four families 
(Witt-Jacobson, Special, Hamiltonian, Contact) of
restricted simple algebras of Cartan-type (see \cite{VIV1} and \cite{VIV2}),
showing that these are never rigid.

By standard facts of deformation theory, the infinitesimal deformations of a
Lie algebra are parametrized by
the second cohomology of the Lie algebra with values in the adjoint representation.
Assuming the notations from sections 2 about the restricted Melikian
algebra $M$ as well as the definition of the squaring operator ${\rm Sq}$, 
we can state the main result of this paper.

\begin{theo}\label{maintheorem}
The infinitesimal deformations of the Melikian algebra $M$ are given by 
$$H^2(M,M)=\langle {\rm Sq}(1)\rangle_F \bigoplus_{i=1}^2 \langle {\rm Sq}(D_i)
\rangle_F \bigoplus_{i=1}^2 \langle {\rm Sq}(\ti{D_i})\rangle_F.$$
\end{theo}

The strategy of our proof consists in exploiting the Hochschild-Serre spectral 
sequence relative to the subalgebra of negative elements of $M$. A similar strategy
has been used by Kuznetsov (\cite{KUT1}) in order to prove that the Melikian
algebra does not admit filtered deformations.
As a byproduct of our proof, we give a new proof (see Theorem \ref{derivations})
of the vanishing of the first cohomology group of the adjoint representation 
(\cite[Prop. 2.2.13]{KUT1}).

The referee suggested another possible approach to prove the above Main Theorem.
Namely, the restricted Melikian algebra $M$ can be realized as a subalgebra of the 
restricted contact algebra $K(5)$ of rank $5$ (see \cite{MEL}) in such a way 
that the negative elements of $M$ coincide with the negative elements of $K(5)$.
This allows to reduce the above result to the analogous result for the contact algebras
(see \cite{VIV2}).

The paper is organized as follows. In Section 2 we recall, in order to fix the 
notations, the basic properties of the restricted Melikian algebra, the
Hochschild-Serre spectral sequence and the squaring operators. In Section 3
we prove that the cocycle appearing in the Main Theorem \ref{maintheorem}
are independent in $H^2(M,M)$ and we outline the strategy to prove that they 
generate the whole second cohomology group. Each of the remaining four sections is devoted 
to carry over one of the four steps of this strategy.

The result presented here constitute part of my doctoral thesis. I thank my advisor
prof. Schoof for useful advices and constant encouragement. I thank the referee
for useful suggestions that helped to improve and clarify the exposition.

\section{Notations}

\subsection{The restricted Melikian algebra $M$}

We fix a field $F$ of characteristic $p=5$.
Let $A(2)=F[x_1,x_2]/(x_1^p,x_2^p)$ be the $F$-algebra
of truncated polynomials in $2$ variables and let $W(2)={\rm Der}_F A(2)$
the restricted Witt-Jacobson Lie algebra of rank $2$.
The Lie algebra $W(2)$ is a free $A(2)$-module 
with basis $D_1:=\frac{\partial}{\partial x_1}$ and $D_2:=\frac{\partial}{\partial x_2}$.
Let $\ti{W(2)}$ be a copy of $W(2)$
and for an element $D\in W(2)$ we indicate with $\widetilde{D}$ the corresponding element
inside $\widetilde{W(2)}$. The Melikian algebra $M$ is defined as
$$M=A(2)\oplus W(2)\oplus \widetilde{W(2)},$$
with Lie bracket defined by the following rules (for all $D, E\in W(2)$ and 
$f, g\in A(2)$):
$$\begin{sis}
&[D,\widetilde{E}]:=\widetilde{[D,E]}+2\,{\rm div}(D)\widetilde{E},\\
&[D,f]:=D(f)-2\, {\rm div}(D)f,\\
&[f_1\widetilde{D_1}+f_2\widetilde{D_2},g_1\widetilde{D_1}+g_2\widetilde{D_2}]:=
f_1g_2-f_2g_1,\\
&[f,\widetilde{E}]:=f E,\\
&[f,g]:=2\,(gD_2(f)-fD_2(g))\widetilde{D_1}+2\,(fD_1(g)-gD_1(f))\widetilde{D_2},\\
\end{sis}$$
where ${\rm div}(f_1D_1+f_2D_2):=D_1(f_1)+D_2(f_2)\in A(2)$.
The Melikian algebra $M$ is a restricted simple Lie algebra of dimension $125$ (see
\cite[section 4.3]{STR}) with a $\Z$-grading given by (for all $D, E\in W(2)$ and 
$f\in A(2)$):
$$\begin{sis}
&{\rm deg}_M(D):=3\,{\rm deg}(D),\\
&{\rm deg}_M(\widetilde{E}):=3\,{\rm deg}(E)+2,\\
&{\rm deg}_M(f):=3\,{\rm deg}(f)-2.\\
\end{sis}$$
The lowest terms of the gradation are
$$M_{-3}=FD_1+FD_2,\hspace{0,5cm} M_{-2}=F\cdot 1, \hspace{0,5cm} M_{-1}=F\widetilde{D_1}+
F\widetilde{D_2}, \hspace{0,5cm} M_{0}=\sum_{i, j=1, 2}Fx_iD_j  $$
while the highest term is
$M_{23}=x_1^4x_2^4(FD_1+FD_2).$
Moreover $M$ has a $\Z/3\, \Z$-grading given by:
$$M_{\overline{1}}:= A(2), \hspace{0,5cm} M_{\overline{0}}:=W(2), \hspace{0,5cm}
M_{\overline{2}}:=\widetilde{W(2)}.$$
The Melikian algebra $M$ has a root space decomposition with respect to a canonical Cartan
subalgebra.
\begin{pro}\label{root-deco}
\begin{itemize}
\item[(a)] $T_M:=\langle x_1D_1\rangle_F\oplus \langle x_2D_2\rangle_F$ is a maximal 
torus of $M$ (called the canonical maximal torus).
\item[(b)] The centralizer of $T_M$ inside $M$ is the subalgebra 
$$C_M=T_M\oplus \langle x_1^2x_2^2 \rangle_F\oplus \langle x_1^4x_2^3\widetilde{D_1}
\rangle_F \oplus \langle x_1^3x_2^4\widetilde{D_2}\rangle_F,$$
which is hence a Cartan subalgebra of $M$ (called the canonical Cartan subalgebra).
\item[(c)] Let $\Phi_M:={\rm Hom}_{\mathbb{F}_5}(\oplus_{i=1}^2\langle x_iD_i\rangle_
{\mathbb{F}_5},\mathbb{F}_5)$ where $\mathbb{F}_5$ is the prime field
of $F$. There is a Cartan decomposition
$M=\bigoplus_{\phi\in \Phi_M}M_{\phi}$, where every summand $M_{\phi}$ has dimension 
$5$ over $F$. Explicitly:
$$\begin{sis}
& x^a\in M_{(a_1-2,a_2-2)},\\
& x^a D_i\in M_{(a_1-\delta_{1i},a_2-\delta_{2i})},\\
&\ti{x^aD_i}\in M_{(a_1+2-\delta_{1i}, a_2+2-\delta_{2i})}.
\end{sis}$$
In particular, if $E\in M_{(\phi_1,\phi_2)}$ then ${\deg E}\equiv 3(\phi_1+\phi_2) \mod 5$.
\end{itemize}
\end{pro}
\begin{proof}
See \cite[section 4.3]{STR}.
\end{proof}

\subsection{Cohomology of Lie algebras}

If $\g$ is a Lie algebra over a field $F$ and $M$ is a $\g$-module, then
the cohomology groups $H^*(\g,M)$ can be computed from the
complex of $n$-dimensional cochains $C^n(\g, M)$ ($n\geq 0$),
that are alternating $n$-linear functions $f:\Lambda^{n}(\g)\to M$,
with differential $d:C^n(\g,M)\to C^{n+1}(\g,M)$ defined by
\begin{equation} \label{d} \begin{split}
\d f(\sigma_0,\dots,\sigma_n)=&\sum_{i=0}^n(-1)^i\sigma_i\cdot f(\sigma_0,\dots, \hat{\sigma_i},
\dots, \sigma_n)+\\
&\sum_{p<q}(-1)^{p+q}f([\sigma_p,\sigma_q],\sigma_0,\dots,\hat{\sigma_p},\dots, \hat{\sigma_q},
\dots \sigma_n)
\end{split} \end{equation}
where the sign $\hat{}$ means that the argument below must be omitted. Given 
$f \in C^n(\g, M)$ and $\gamma \in \g$, we denote with $f_{\gamma}$ the restriction
of $f$ to $\gamma\in \g$, that is the element of $C^{n-1}(\g, M)$ given by 
$$f_{\gamma}(\sigma_0, \cdots, \sigma_{n-1}):=f(\gamma, \sigma_0, \cdots, 
\sigma_{n-1}).$$
With this notation, the above differential satisfies the 
following useful formula (for any $\gamma\in \g$ and $f\in C^n(\g,M)$):
\begin{gather}
\d(\gamma \cdot f)=\gamma\cdot (\d f)\label{module-d}, \\
(\d f)_{\gamma}=\gamma\cdot f-\d(f_{\gamma}),\label{restriction}
\end{gather}
where each $C^n(\g,M)$ is a $\g$-module by means of the action
\begin{equation}\label{module}
(\gamma\cdot f)(\sigma_1,\dots,\sigma_n)=\gamma\cdot f(\sigma_1,\dots, \sigma_n)
-\sum_{i=1}^n f(\sigma_1,\cdots,[\gamma,\sigma_i],\dots \sigma_n).
\end{equation}
As usual we indicate with $Z^n(\g,M)$ the subspace of $n$-cocycles and with
$B^n(\g,M)$ the subspace of $n$-coboundaries. Therefore
$H^n(\g,M):=Z^n(\g,M)/B^n(\g,M)$.

A useful tool to compute cohomology of Lie algebras is the 
Hochschild-Serre spectral sequence relative to a subalgebra.
Given a subalgebra $\h < \g$, one can define a decreasing filtration
$\{F^j C^n(\g, M)\}_{j=0, \cdots, n+1}$ on the space of $n$-cochains:
$$
F^j C^n(\g, M)=\{f\in C^n(\g, M)\: | \: f(\sigma_1, \cdots, \sigma_n)=0\:
\text{ if } \sigma_1,\cdots, \sigma_{n+1-j}\in \h \}.
$$
This gives rise to a spectral sequence converging to the cohomology $H^n(\g, M)$,
whose first level is equal to (see \cite{HS}):
\begin{equation}\label{HS-ss1}
E_1^{p,q}=H^q(\h, C^p(\g/\h, M))\Longrightarrow H^{p+q}(\g, M).
\end{equation}
In the case where $\h$ is an ideal of $\g$ (which we indicate as $\h\lhd \g$) 
the above spectral sequence becomes
\begin{equation}\label{HS-ss2}
E_2^{p,q}=H^p(\g/\h, H^q(\h, M))\Longrightarrow H^{p+q}(\g, M).
\end{equation}
Moreover for the second page of the first spectral sequence
(\ref{HS-ss1}), we have the equality
\begin{equation}\label{REL-coho}
E_2^{p,0}= H^p(\g,\h;M),
\end{equation}
where $H^*(\g,\h;M)$ are the relative cohomology groups defined (by
Chevalley and Eilenberg \cite{CE}) from the sub-complex
$C^p(\g,\h;M)\subset C^p(\g,M)$ consisting of cochains orthogonal to
$\h$, that is cochains satisfying the two conditions:
\begin{gather}
f_{| \h}=0, \label{ortho1}\\
\d f_{| \h}=0 \hspace{0,5cm} \text{ or equivalently } \hspace{0,5cm}
\gamma \cdot f=0 \hspace{0,3cm} \text{ for every }
\gamma\in \h.\label{ortho2}
\end{gather}
Note that in the case where $\h\lhd \g$, the equality (\ref{REL-coho})
is consistent with the second spectral sequence (\ref{HS-ss2}) because in that case
we have $H^p(\g, \h, M)=H^p(\g/\h, M^{\h})$.

Suppose that a torus $T$ acts on both $\g$ and $M$ in a way that is
compatible with the action of $\g$ on $M$, which means that $t\cdot(g\cdot m)=(t\cdot g)\cdot m+
t\cdot (g\cdot m)$ for every $t\in T$, $g\in \g$ and $m\in M$. Then the
action of $T$ can be extended to the space of $n$-cochains by
$$(t\cdot f)(\sigma_1,\cdots ,\sigma_n)=t\cdot f(\sigma_1, \cdots, \sigma_n)
-\sum_{i=1}^n f(\sigma_1, \cdots, t\cdot \sigma_i, \cdots, \sigma_n).$$
It follows easily from the compatibility of the action of $T$ and formula
\ref{restriction}, that the action of $T$ on the cochains commutes with the
differential $\d$. Therefore, since the action of a torus is always completely
reducible, we get a decomposition in eigenspaces
\begin{equation}\label{deco1}
H^n(\g,M)=\bigoplus_{\phi\in \Phi} H^n(\g,M)_{\phi}
\end{equation}
where $\Phi={\rm Hom}_F(T,F)$ and $H^n(\g,M)_{\phi}=\{[f] \in H^n(\g,M)\: | \:
t\cdot [f]=\phi(t) [f] \: \text{ if }  t\in T\}$. A particular case of this situation
occurs when $T\subset \g$ and $T$ acts on $\g$ via the adjoint action
and on $M$ via restriction of the action of $\g$. It is clear that
this action is compatible and moreover the above decomposition reduces to
$$H^n(\g,M)=H^n(\g,M)_{\0}$$
where $\underline{0}$ is the trivial homomorphism (in this situation we say that
the cohomology reduces to \emph{homogeneous} cohomology). Indeed, if we consider
an element $f\in Z^n(\g,M)_{\phi}$, then by applying formula \ref{restriction}
with $\gamma=t\in T$ we get
$$0=(\d f)_t=t\cdot f- \d (f_t)=\phi(t)f-\d(f_t),$$
from which we see that the existence of a $t\in T$ such that $\phi(t)\neq 0$
forces $f$ to be a coboundary.

Now suppose that $\g$ and $M$ are graded and that the action of $\g$ respects
these gradings, which means that $\g_d\cdot M_e\subset M_{d+e}$ for all $e, d\geq 0$.
Then the space
of cochains can also be graded: a homogeneous cochain $f$ of degree $d$ is
a cochain such that $f(\g_{e_1}\times \cdots \times \g_{e_n})\subset M_{\sum e_i+d}$.
With this definition, the differential becomes of degree $0$ and therefore we get
a degree decomposition
\begin{equation}\label{deco2}
H^n(\g,M)=\bigoplus_{d\in \Z} H^n(\g,M)_d.
\end{equation}
Finally, if the action of $T$ is compatible with the grading, in the sense
that $T$ acts via degree $0$ operators both on $\g$ and on $M$, then the above
two decompositions \ref{deco1} and \ref{deco2} are compatible and give rise to the refined
weight-degree decomposition
\begin{equation}\label{deco3}
H^n(\g,M)=\bigoplus_{\phi\in \Phi}\bigoplus_{d\in \Z}H^n(\g,M)_{\phi,d}.
\end{equation}

We will use frequently the above weight-degree decomposition with respect to the 
action of the maximal torus $T_M\subset M_0$ of the restricted Melikian 
algebra $M$ (see Proposition \ref{root-deco}).

\subsection{Squaring operation}

There is a canonical way to produce $2$-cocycles in $Z^2(\g,\g)$ over a field of characteristic
$p>0$, namely the squaring operation (see \cite{GER1}).
Given a derivation $\gamma\in Z^1(\g,\g)$ (inner or not), one defines the squaring of
$\gamma$ to be
\begin{equation}\label{Square}
\Sq(\gamma)(x,y)=\sum_{i=1}^{p-1}\frac{[\gamma^i(x),\gamma^{p-i}(y)]}{i!(p-i)!}\in Z^2(\g,\g)
\end{equation}
where $\gamma^i$ is the $i$-iteration of $\gamma$.
In \cite{GER1} it is shown that $[\Sq(\gamma)]\in H^2(\g,\g)$ is an obstruction to integrability
of the derivation $\gamma$, that is to the possibility of finding an automorphism
of $\g$ extending the infinitesimal automorphism given by $\gamma$.

\section{Strategy of the proof of the Main Theorem}

First of all, we show that the five cocycles $\{{\rm Sq}(\ti{D_i})\}_{i=1, 2}$,
$\{{\rm Sq}(1)\}$, $\{{\rm Sq}(D_i)\}_{i=1, 2}$ are independent in $H^2(M,M)$.
Observe that the first two cocycles have degree $-5$, the third has degree $-10$ 
and the last two have degree $-15$. Therefore, according to the decomposition 
(\ref{deco2}), it is enough to show that the first two are independent, the third 
is non-zero and the last two are independent in $H^2(M, M)$.  

To prove the independence of the first two cocycles, we observe that (for $i\neq j$) 
\begin{equation}\label{indep1}
{\rm Sq}(\ti{D_i})(x_iD_j, x_i \ti{D_j})=D_i,
\end{equation}
while for a cochain $g\in C^1(M,M)_{\0, -5}$ we have that
\begin{equation}\label{indep1bis}
\d g(x_iD_j, x_i \ti{D_j})=[x_iD_j, g(x_i\ti{D_j})]-[x_i\ti{D_j},g(x_i D_j)]-
g([x_iD_j, x_i \ti{D_j}])=0,
\end{equation}
since $g(x_iD_j)=0$ by degree reasons, $g(x_i\ti{D_j})=0$ by degree and homogeneity 
reasons and $[x_iD_j, x_i \ti{D_j}]=0$.  

The third cocyle is non-zero because
\begin{equation}\label{indep2}
\begin{sis}
&{\rm Sq}(1)(x_i, x_ix_j^2D_i)=2D_i,\\
&{\rm Sq}(1)(x_i, x_j^3D_j)=D_i,\\
\end{sis}
\end{equation}
while for a cochain $g\in C^1(M, M)_{\0, -10}$ we have that
\begin{equation}\label{indep2bis}
\begin{sis}
&\d g(x_i, x_ix_j^2D_i)=-g(x_ix_j^2),\\
&\d g(x_i, x_j^3D_j)=-g(x_ix_j^2).\\
\end{sis}
\end{equation}

Finally, the independence of the last two cocycles follows from 
\begin{equation}\label{indep3}
\begin{sis}
&{\rm Sq}(D_i)(x_i D_j, x_i^4x_jD_j)=D_j,\\
&{\rm Sq}(D_i)(x_i D_j, x_i^3x_j^2D_j)=0,\\
\end{sis}
\end{equation}
together with the fact that for a cochain $g\in C^1(M, M)_{\0, -15}$ we have that
\begin{equation}\label{indep3bis}
\begin{sis}
&\d g(x_iD_j, x_i^4x_j^2D_i)=[x_iD_j, g(x_i^4x_jD_j)],\\
&\d g(x_iD_j, x_i^3x_j^2D_j)=-g([x_iD_j, x_i^3x_j^2D_j])=-2g(x_i^4x_jD_j).\\
\end{sis}
\end{equation}

The remaining of this paper is devoted to show that the above five cocycles generate
the cohomology group $H^2(M,M)$, as stated in the Main Theorem \ref{maintheorem}. 
We outline here the strategy of the proof that will be carried over in the next 
sections. The proof is divided in five steps: 

\underline{STEP I}: We prove in Corollary \ref{M-firststep} that 
$$H^2(M,M)=H^2(M,M_{<0};M).$$

\underline{STEP II}: We prove in Proposition \ref{M-rel-coho} that 
$$H^2(M,M_{<0};M)\subset H^2(M_{\geq 0},M_{-3}),$$
where $M_{\geq 0}$ acts on $M_{-3}=\langle D_1,D_2\rangle_F$ via projection onto
$M_{\geq 0}/M_{\geq 1}=M_0$ followed by the adjoint representation of $M_0$
onto $M_{-3}$.

\underline{STEP III}: We prove in Corollary \ref{reduction>0} that 
$$H^2(M_{\geq 0},M_{-3})\subset
\bigoplus_{i=1}^2 \langle {\rm Sq}(\ti{D_i})\rangle_F\bigoplus 
H^2(M_{\geq 1},M_{-3})^{M_0},$$
where $M_{\geq 1}$ acts trivially on $M_{-3}$.

\underline{STEP IV}: We prove in Proposition \ref{inva-coho} that 
$$H^2(M_{\geq 1},M_{-3})^{M_0}= \bigoplus_{i=1}^2 \langle 
{\rm Sq}(D_i)\rangle_F \bigoplus \langle {\rm Sq}(1)\rangle_F.$$

As a byproduct of our Main Theorem, we obtain a new proof of the following result
(\cite[Prop. 2.2.13]{KUT1}, see also \cite[chapter 7]{STR}).

\begin{theo}[Kuznetsov]\label{derivations}
$H^1(M,M)=0$.
\end{theo} 
\begin{proof}
The spectral sequence (\ref{HS-negative}), together with Proposition \ref{M_{<0}-cohomology},
gives that $$H^1(M,M)=E_2^{1,0}=H^1(M,M_{<0};M).$$
The Proposition \ref{M-rel-coho} gives that
$H^1(M,M_{<0};M)\subset H^1(M_{\geq 0},$ $M_{-3}).$
Using the spectral sequence (\ref{HS-0term}), together with Propositions \ref{M0-level} 
and \ref{M1-level}, we get the vanishing of this last group.
\end{proof}

\section{Step I: Reduction to $M_{<0}$-relative cohomology}

In this section, we carry over the first step outlined in section 3.
To this aim, we consider 
the homogeneous Hochschild-Serre spectral sequence associated 
to the subalgebra $M_{<0}<M$ (see section 2.2):
\begin{equation}\label{HS-negative}
(E_1^{r,s})_{\0}=H^s(M_{<0},C^r(M/M_{<0},M))_{\0}\Rightarrow H^{r+s}(M, M)_{\0}=
H^{r+s}(M,M).
\end{equation}
We adopt the following notation: given elements
$E_1,\cdots, E_n\in M_{<0}$ and $G\in M$ we denote with $\delta_{E_1,\cdots,E_n}^G$ the
cochain of $C^n(M_{<0},M)$ whose only non-zero values are
$$\delta_{E_1, \cdots, E_n}^G(E_{\sigma(1)},\cdots,E_{\sigma(n)})={\rm sgn}(\sigma) G,
$$  
for any permutation $\sigma\in S_n$.

\begin{pro}\label{M_{<0}-cohomology}
In the above spectral sequence (\ref{HS-negative}), we have 
$$(E_1^{0,1})_{\0}=(E_1^{0,2})_{\0}=0.$$
\end{pro}
\begin{proof}
Consider the Hochschild-Serre spectral sequence associated 
to the ideal $M_{-3}\lhd M_{\leq -2}$:
\begin{equation}\label{HS2}
E_2^{r,s}=H^r(M_{\leq -2}/M_{-3},H^s(M_{-3},M))\Rightarrow H^{r+s}(M_{\leq -2},M).
\end{equation}
The lowest term
$M_{-3}=\langle D_1, D_2\rangle_F$ acts on $M=A(2)\oplus 
W(2)\oplus \widetilde{W(2)}$
via its natural action on $A(2)$ and via adjoint representation on $W(2)$ and 
$\widetilde{W(2)}$.
Hence, according to \cite[Prop. 3.4 and Cor. 3.5]{VIV1}, we have that (for $s=0, 1, 2$)
\begin{equation*}
H^s(M_{-3},M)=
\begin{sis}
& M_{<0} & \text{ if } s=0,\\
& \bigoplus_{G\in M_{<0}}\langle \delta_{D_1}^{x_1^4 G}, \delta_{D_2}^{x_2^4 G} 
\rangle_F & \text{ if } s=1,\\
& \bigoplus_{G\in M_{<0}}\langle \delta_{D_1, D_2}^{x_1^4x_2^4 G} \rangle_F & 
\text{ if } s=2.\\
\end{sis}
\end{equation*}
Moreover $M_{\leq -2}/M_{-3}=\langle 1 \rangle_F$ acts on the above cohomology groups
$H^s(M_3,M)$ via its adjoint action on 
$M_{<0}=\langle 1 \rangle \oplus \langle D_1, D_2 \rangle \oplus \langle \widetilde{D_1},
\widetilde{D_2} \rangle$, that is
$$[1,1]=0, \hspace{0,5cm} [1,D_i]=0, \hspace{0,5cm}
[1,\widetilde{D_i}]=D_i.$$
Hence, using that the above Hochschild-Serre spectral sequence (\ref{HS2}) is 
degenerate since $M_{\leq -2}/M_{-3}=\langle 1\rangle_F$ has dimension $1$,
we deduce that
$$H^s(M_{\leq -2},M)=
\begin{sis}
& M_{-3}\oplus M_{-2} & \text{ if } s=0,\\
& \bigoplus_{G\in M_{-3}\oplus M_{-2}}\langle \delta_{D_1}^{x_1^4 G}, 
\delta_{D_2}^{x_2^4 G} 
\rangle_F \bigoplus_{H\in M_{-2}\oplus M_{-1}} \langle \delta_1^H \rangle_F & \text{ if } s=1,\\
&  \bigoplus_{G\in M_{-3}\oplus M_{-2}}\langle \delta_{D_1, D_2}^{x_1^4 x_2^4 G} 
\rangle_F \bigoplus_{H\in M_{-2}\oplus M_{-1}} \langle \delta_{1, D_1}^{x_1^4 H}, 
\delta_{1, D_2}^{x_2^4 H} \rangle_F& \text{ if } s=2.\\
\end{sis}$$

Finally consider the homogeneous Hochschild-Serre spectral sequence associated to the ideal
$M_{\leq -2}\lhd M_{<0}$:
\begin{equation}\label{HS3}
(E_2^{r,s})_{\0}=H^r(M_{<0}/M_{\leq -2},H^s(M_{\leq -2},M))_{\0} \Rightarrow
H^{r+s}(M_{<0},M)_{\0}.
\end{equation}
From the explicit description of above, one can easily check that the only
non-zero terms and non-zero maps of the above spectral sequence that can contribute 
to $H^s(M_{<0},M)_{\0}$ for $s=1,2$ are
\begin{align*}
(E_2^{0,1})_{\0}=\langle \delta_1^1\rangle_F &\rightarrow (E_2^{2,0})_{\0}=
\langle \delta_{\widetilde{D_1}, \widetilde{D_2}}^1\rangle_F, \\
(E_2^{0,2})_{\0}=\bigoplus_{i=1}^2 \langle \delta_{D_i, 1}^{x_i^4}\rangle_F
&\rightarrow (E_2^{2,1})_{\0}=\bigoplus_{i=1}^2 \langle \delta_{D_i,
\widetilde{D_1},\widetilde{D_2}}^{x_i^4}\rangle_F, 
\end{align*}
where the maps are given by the differentials.
Using the relation $[\widetilde{D_1},\widetilde{D_2}]=1$, it is easy to see that all the 
above maps are isomorphisms and hence the conclusion follows.
\end{proof}

In the next Proposition, we need the following 

\begin{lem}\label{M-nonhom}
We have that
$$H^1(M_{<0},M)=\bigoplus_{k, h}\langle \delta_{D_k}^{x_k^4D_h}\rangle_F
\bigoplus_{i\neq j}\langle \delta_{\ti{D_i}}^{D_j}\rangle_F \bigoplus
\langle \delta_{\ti{D_1}}^{D_1}-\delta_{\ti{D_2}}^{D_2}\rangle_F.
$$
\end{lem}
\begin{proof}
From the (non-homogeneous) Hochschild-Serre spectral sequence  associated to the ideal
$M_{\leq -2}\lhd M_{<0}$:
\begin{equation*}
(E_2^{r,s})=H^r(M_{<0}/M_{\leq -2},H^s(M_{\leq -2},M)) \Rightarrow
H^{r+s}(M_{<0},M),
\end{equation*}
we deduce the exact sequence
$$0\to H^1(M_{<0}/M_{\leq -2},M^{M_{\leq -2}})\to H^1(M_{<0},M)\to H^1(M_{\leq -2},M)^
{M_{<0}/M_{\leq -2}}\stackrel{\partial}{\longrightarrow} $$
$$\stackrel{\partial}{\longrightarrow} H^2(M_{<0}/M_{\leq -2},M^{M_{\leq -2}}).$$
Using the computation of $H^s(M_{\leq -2}, M)$ for $s=0, 1$ from the Proposition
\ref{M_{<0}-cohomology}, it is easily seen that
$$\begin{sis}
& H^1(M_{<0}/M_{\leq -2},M^{M_{\leq -2}})=\frac{\bigoplus_{i,j}
\langle \delta_{\widetilde{D_i}}^{D_j}\rangle_F }{\langle {\rm ad}(1)\rangle_F},\\
& H^1(M_{\leq -2},M)^{M_{<0}/M_{\leq -2}}=\bigoplus_{h,k}\langle \delta_{D_h}^{x_h^4 D_k}\rangle_F
\oplus \langle \delta_1^1 \rangle_F,\\
& H^2(M_{<0}/M_{\leq -2},M^{M_{\leq -2}})=\langle \delta_{(\widetilde{D_1},\widetilde{D_2})}^1
\rangle_F.\\
\end{sis}$$
Using the relation $[\widetilde{D_1},\widetilde{D_2}]=1$, it's easy to see
that $\partial(\delta_1^1)=\delta_{(\widetilde{D_1},\widetilde{D_2})}^1$ which gives the
conclusion.
\end{proof}

\begin{pro}\label{ME_2^{1,1}}
In the above spectral sequence (\ref{HS-negative}), we have that
$(E_2^{1,1})_{\0}=0.$
\end{pro}
\begin{proof}
We have to show the injectivity of the level $1$ differential map
$$d:(E_1^{1,1})_{\0}\longrightarrow(E_1^{2,1})_{\0}.$$
In the course of this proof, we adopt the following convention: given an element
$f\in  C^1(M_{<0}, C^s(M/M_{<0},M)$,
 we write its value on $D\in M_{<0}$ as $f_{D}\in C^s(M/M_{<0},M)$.

We want to show, by induction on the degree of $E\in M/M_{<0}$,
that if $[df]=0\in H^1(M_{<0},C^2(M/M_{<0},M))$ then we can choose a representative
$\widetilde{f}$ of $[f]\in H^1(M_{<0},C^1(M/M_{<0},M))$ such that
$\widetilde{f}_{D}(E)=0$ for every $D\in M_{<0}$. So suppose that we have already found
a representative $f$
 such that $f_{D}(F)=0$ for every $F\in M/M_{<0}$ of degree less than $d$
and for every $D\in M_{<0}$.
First of all, we can find a representative $\widetilde{f}$ of $[f]$ such that
\begin{equation*}\tag{*}
\begin{sis}
& \widetilde{f}_{D_i}(E)=\sum_{j=1}^2\alpha_i^j x_i^{4}D_j \text{ for } i=1, 2, \\
& \ti{f}_1(E)=0,\\
& \widetilde{f}_{\ti{D_1}}(E)=\beta_1 D_2+\gamma D_1,\\
& \widetilde{f}_{\ti{D_2}}(E)=\beta_2 D_1-\gamma D_2,\\
\end{sis}
\end{equation*}
for every $E\in M$ of degree $d$. Indeed, by the induction
hypothesis, the cocycle condition for $f$ is $\partial f_{D,D'}(E)=[D,f_{D'}(E)]-
[D',f_{D}(E)]-f_{[D, D']}(E)$ for any $D, D'\in M_{<0}$.
On the other hand, by choosing an element $h\in C^1(M/M_{<0},M)$ that vanishes on
the elements of degree less than $d$, we can add to $f$ (without changing its  cohomological
class neither affecting the inductive assumption) the coboundary $\partial h$ whose value
on $E$ are $\partial h_{D}(E)=[D,h(E)]$.  Hence, for a fixed element $E$ of degree $d$,
the map $D\mapsto f_{D}(E)$ gives rise to an element of $H^1(M_{<0},M)$ and,
by Lemma \ref{M-nonhom}, we can chose an element $h(E)$ as above
such that the new cochain $\widetilde{f}=f+\partial h$ verifies the condition  $(*)$ of above.

Note that, by the homogeneity of $\widetilde{f}$, we have the following pairwise 
disjoint possibilities for $E$:
$$\begin{sis}
& \alpha_i^j\neq 0 \text{ for } i=1, 2 &\Rightarrow E\in \{x_j x_{j+1}^2, 
x_j^4x_{j+1}D_{j+1}, \ti{x_j^3x_{j+1}^3D_j}, \ti{x_j^2x_{j+1}^4D_{j+1}}\},\\
& \beta_j\neq 0 &\Rightarrow E\in \{x_j x_{j+1}^4, x_{j+1}^2D_j,
x_j^4x_{j+1}^3D_{j+1}, \ti{x_j^3 D_j}, \ti{x_j^2x_{j+1}D_{j+1}}\},\\
& \gamma\neq 0 &\Rightarrow E\in \{x_1^4x_{2}^3D_1,
x_1^3x_{2}^4D_{2}, \ti{x_1^2x_2 D_1}, \ti{x_1 x_2^2 D_{2}}\}.
\end{sis}$$

Now we are going to use the condition that $[d\widetilde{f}]=0\in (E_1^{2,1})_{\0}$,
that is $d\widetilde{f}=\partial g$ for some $g\in C^2(M/M_{<0},M)_{\0}$.
Explicitly, for $A, B\in M/M_{<0}$ we have that (for any $D\in M_{<0}$)
\begin{equation}\label{Wfor-g}
\partial g_{D}(A,B)=[D,g(A,B)]-g([D,A],B)-g(A,[D,B]),
\end{equation}
\begin{equation}\label{Wfor-d}
d\widetilde{f}_{D}(A,B)=
\widetilde{f}_{D}([A,B])-[A,\widetilde{f}_{D}(B)]+[B,\widetilde{f}_{D}(A)]-
\widetilde{f}_{[D,A]}(B)+\widetilde{f}_{[D,B]}(A),
\end{equation}
where the last two terms in the first formula can be non-zero only if $[D, A]\in M_{\geq 0}$ and
$[D, B]\in M_{\geq 0}$ respectively,
where the last two terms in the second formula can be non-zero only if $[D, A]\in M_{<0}$ and
$[D, B]\in M_{<0}$ respectively.

Suppose first that $\alpha_i^j\neq 0$ for a fixed $j$ and for $i=1, 2$. Then it is 
straightforward to check that for each $E$ in the above list it is possible to find 
two elements $A, B\in M$ such that $[A, B]=E$, ${\rm deg}(B)=0$ and $A$ does not belong 
to the above list.
Apply the above formulas for each such pair $(A, B)$.
Taking into account the inductive hypothesis on the degree and the homogeneity assumptions,
the formula (\ref{Wfor-d}) becomes
$$d\widetilde{f}_{D_i}(A,B)=\widetilde{f}_{D_i}([A, B])=\widetilde{f}_{D_i}(E)=\alpha_i^j
x_i^4D_j,$$
while the formula (\ref{Wfor-g}) gives
$$\partial g_{D_i}(A,B)=[D_i,g(A,B)]-g([D_i,A],B).$$
The first term in the last expression is a derivation with respect to $D_i$ and therefore
it cannot involve the monomial $x_i^4D_j$. The same is true for the second term 
as it follows by applying the formulas (\ref{Wfor-d}) and (\ref{Wfor-g}) for 
the elements ${\rm ad}(D_i)^k(A)$ (with $k=1, \cdots, p-1$) 
and $B$ and using the vanishing hypothesis:
$$\begin{sis}
&d \widetilde{f}_{D_i}({\rm ad}(D_i)^k(A),B)=0,\\
&\partial g_{D_i}({\rm ad}(D_i)^k(A), B)=[D_i, g({\rm ad}(D_i)^k(A),B)]-
g({\rm ad}(D_i)^{k+1}(A), B).
\end{sis}$$
Therefore we conclude that $\alpha_i^j=0$, an absurd. The other cases $\beta_j\neq 0$ and 
$\gamma\neq 0$ are excluded using a similar argument.
\end{proof}

Finally we get the main result of this section.

\begin{cor}\label{M-firststep}
We have that  
$$H^2(M,M)=H^2(M,M_{<0};M).$$
\end{cor}
\begin{proof}
From the spectral sequence (\ref{HS-negative}), using the vanishing of $(E_1^{0,2})_{\0}$
(Proposition \ref{M_{<0}-cohomology}) and of $(E_1^{1,1})_{\0}$ (Proposition 
\ref{ME_2^{1,1}}), we get that 
$$H^2(M,M)=(E_{\infty}^{2,0})_{\0}=(E_2^{2,0})_{\0}=H^2(M,M_{<0};M).$$
\end{proof}

\section{Step II: Reduction to $M_{\geq 0}$-cohomology}

In this section, we carry over the second step of proof of the main theorem
(see section 3). 

Consider the action of $M_{\geq 0}$ on $M_{-3}=\langle D_1,D_2\rangle_F$ obtained 
via projection onto
$M_{\geq 0}/M_{\geq 1}=M_0$ followed by the adjoint representation of $M_0$
onto $M_{-3}$.

\begin{pro}\label{M-rel-coho}
We have that 
$$\begin{sis}
& H^1(M, M_{<0}; M)\subset H^1(M_{\geq 0}, M_{-3}), \\
& H^2(M, M_{<0}; M)\subset H^2(M_{\geq 0}, M_{-3}).
\end{sis}$$
\end{pro}
\begin{proof}
For every $s\in \Z_{\geq 0}$, consider the map
$$\phi_s: C^s(M,M_{< 0};M)
\to C^s(M_{\geq 0},M_{-3})$$ induced by the restriction to the
subalgebra $M_{\geq 0}\subset M$ and by the projection
$M\twoheadrightarrow M/M_{\geq -2}=M_{-3}$. It is
straightforward to check that the maps $\phi_s$ commute with the
differentials and hence they define a map of complexes. Moreover the
orthogonality conditions with respect to the subalgebra $M_{<0}$
give the injectivity of the maps $\phi_s$. Indeed, on one hand, the
condition (\ref{ortho1}) says that an element $f\in
C^s(M,M_{<0};M)$ is determined by its restriction to
$\wedge^s M_{\geq 0}$. On the other hand, the condition
(\ref{ortho2}) implies that the values of $f$ on an $s$-tuple are
determined, up to elements of $M^{M_{<0}}=M_{-3}$, by
induction on the total degree of the $s$-tuple. 

\noindent The map $\phi_0$ is an isomorphism since  
$$C^0(M, M_{<0}; M)=M^{M_{<0}}=M_{-3}=C^0(M_{\geq 0}, M_{-3}).$$
From this we get the first statement of the Proposition.
Moreover, it is easily checked that if $g\in C^1(M_{\geq 0},M_{-3})$ is such that
$\d g\in C^2(M,M_{<0};M)$, then $g\in C^1(M, M_{<0}; M)$. This gives
the second statement of the Proposition.
\end{proof}

\section{Step III: Reduction to $M_0$-invariant cohomology}

In this section, we carry over the third step of the proof of the Main Theorem
(see section 3). To this aim, we consider the Hochschild-Serre spectral sequence
relative to the ideal $M_{\geq 1}\lhd M_{\geq 0}$:
\begin{equation}\label{HS-0term}
E_2^{r,s}=H^r(M_0,H^s(M_{\geq 1},M_{-3}))\Longrightarrow H^{r+s}(M_{\geq 0},M_{-3}).
\end{equation}
The first line $E_2^{*,0}$ of the above spectral sequence vanish.  

\begin{pro}\label{M0-level}
In the above spectral sequence (\ref{HS-0term}), we have for every $r\geq 0$:
$$E_2^{r,0}=H^r(M_0,M_{-3})=0.$$
\end{pro}
\begin{proof}
Observe that, since the canonical maximal torus $T_M$ is contained in $M_0$, 
we can restrict to homogeneous cohomology (see section 2.2). But there are 
no homogeneous cochains in $C^r(M_0,M_{-3})$. 
Indeed the weights that occur in $M_{-3}$ are $-\epsilon_i$ while the
weights that occur in $M_0$ are $\0$ and $\epsilon_i-\epsilon_j$.
Therefore the weights that occur in $M_0^{\otimes k}$ have degree
congruent to $0$ modulo $5$ and hence they cannot be equal to
$-\epsilon_i$.
\end{proof}

Next, we determine the groups $E_2^{r,1}=H^r(M_0,H^1(M_{\geq 1},M_{-3}))$ for 
$r=1, 2$ of the above spectral sequence (\ref{HS-0term}).

\begin{pro}\label{M1-level}
In the above spectral sequence (\ref{HS-0term}), we have that
$$E_2^{r,1}=
\begin{sis}
& 0 & \text{ if }  r=0,\\
& \langle \overline{{\rm Sq}(\widetilde{D_1})}, \overline{{\rm Sq}(\widetilde{D_2})}\rangle_F
  & \text{if } r=1,\\
\end{sis}$$
where $\overline{{\rm Sq}(\widetilde{D_i})}$ (for $i=1,2$) is the restriction
of ${\rm Sq}(\widetilde{D_i})$ to $M_0\times M_2$.
\end{pro}
\begin{proof}
Using the Lemma \ref{M-commutators} below, we have that
$$H^1(M_{\geq 1},M_{-3})=\{f\in C^1(M_1+M_2,M_{-3})\: | \: f(\widetilde{x_1D_1})=
-f(\widetilde{x_2D_2})\}.$$

Observe that, since the canonical maximal torus $T_M$ is contained in $M_0$, 
we can restrict to homogeneous cohomology (see section 2.2).
The vanishing of $E_2^{0,1}$ follows directly from the fact that there are no homogeneous
cochains in $C^1(M_0,H^1(M_{\geq 1},$ $  M_{-3}))$.  

The elements $\overline{{\rm Sq}(\ti{D_i})}$ for $i=1,2$ belong
to $H^1(M_0,H^1(M_{\geq 1},M_{-3}))$ and are non-zero 
in virtue of the formulas (\ref{indep1}) and (\ref{indep1bis}).
Moreover it is easy to see that, for
homogeneity reasons, $C^1(M_0,H^1(M_{\geq 1},M_{-3}))_{\0}$ is generated by 
$\overline{{\rm Sq}(\ti{D_i})}$ ($i=1,2$), which gives the conclusion.
\end{proof}

\begin{lem}\label{M-commutators}
$[M_{\geq 1},M_{\geq 1}]=M_{\geq 3}+\langle \widetilde{x_1D_1}+\widetilde{x_2D_2}\rangle_F.$
\end{lem}
\begin{proof}
Clearly $[M_{\geq 1},M_{\geq 1}]\subset M_{\geq 2}$ and
$[M_{\geq 1},M_{\geq 1}]\cap M_2=[M_1,M_1]$.  Since $M_1=\langle x_1, x_2\rangle_F$,
we have that $[M_1, M_1]=\langle [x_1,x_2]\rangle_F=\langle -2(\widetilde{x_1D_1}+
\widetilde{x_2D_2})\rangle_F$.
Hence the proof is complete if we show that
$M_{\geq 3}\subset  [M_{\geq 1},M_{\geq 1}]$. We will consider the $\Z/3\Z$-grading on $M$
and we will consider separately $M_{\geq 3}\cap M_{\overline{i}}$, with $i=0, 1, 2$.
\begin{itemize}

\item[(i)] $M_{\geq 3}\cap M_{\overline{0}}\subset [M_{\geq 1},M_{\geq 1}]$ because 
of the formula $[x_j,\widetilde{x^aD_i}]=x_jx^aD_i.$

\item[(ii)] $M_{\geq 3}\cap M_{\overline{1}}=M_{\geq 4}\cap M_{\overline{1}}\subset
[M_{\geq 1}, M_{\geq 1}].$\\
Indeed, the formula $[x^aD_i,x_i]=(1-2a_i)x^a$ shows that $x^a$ (with $|a|\geq 2$)
belongs to $[M_{\geq 1}, M_{\geq 1}]$ provided that $x^a\neq x_1^3x_2^3$ and this
exceptional case is handled with the formula $[x_1^2x_2^4D_2,x_1]=-2\cdot 4 x_1^3x_2^3$.

\item[(iii)] $M_{\geq 3}\cap M_{\overline{2}}=M_{\geq 5}\cap M_{\overline{2}}\subset
[M_{\geq 1}, M_{\geq 1}].$\\
Indeed, the formula $[x_i^2D_i,\widetilde{x^aD_j}]=\widetilde{[x_i^2D_i,x^aD_j]}+
4\widetilde{x_ix^aD_j}=(a_i+4-\delta_{ij})\widetilde{x^{a+\epsilon_i}D_j}$ shows
the inclusion for the elements of $M_{\geq 5}\cap M_{\overline{2}}$ with the exception
of the two elements of the form $\widetilde{x_i^2x_j^4D_j}$ for $i\neq j$, for which we use
the formula $[x_i^2D_j,\widetilde{x_ix_j^4D_i}]=4\widetilde{x_i^3x_j^3D_i}-
2\widetilde{x_i^2x_j^4D_j}$.
\end{itemize}
\end{proof}

Finally we get the result we were interested in.

\begin{cor}\label{reduction>0}
We have that
$$ H^2(M_{\geq 0},M_{-3})\subset  \bigoplus_{i=1}^2 \langle 
{\rm Sq}(\widetilde{D_i}) \rangle_F \bigoplus H^3(M_{\geq 1}, M_{-3})^{M_0}.$$
\end{cor}
\begin{proof}
From the above spectral sequence (\ref{HS-0term}), using the Propositions \ref{M0-level} and 
\ref{M1-level}, we get the exact sequence
$$0\to \langle {\rm Sq}(\widetilde{D_1}), {\rm Sq}(\widetilde{D_2})
 \rangle_F \to H^2(M_{\geq 0},M_{-3}) \to  H^2(M_{\geq 1},M_{-3})^{M_0}, $$
which gives the conclusion.
\end{proof}

\section{Step IV: Computation of $M_0$-invariant cohomology}

In this section, we carry over the fourth and last step of the proof of the Main Theorem
by proving the following 

\begin{pro}\label{inva-coho}
We have that 
$$H^2(M_{\geq 1},M_{-3})^{M_0}=\bigoplus_{i=1}^2 \langle 
{\rm Sq}(D_i)\rangle_F \bigoplus \langle {\rm Sq}(1)\rangle_F .$$
\end{pro}

\begin{proof}
The strategy of the proof is
to compute, step by step as $d$ increases, the truncated invariant cohomology groups
$$H^2\left(\frac{M_{\geq 1}}{M_{\geq d+1}},M_{-3}\right)^{M_0}.$$  
Observe that if $d\geq 23$ then $M_{\geq d+1}=0$ and
hence we get the cohomology we are interested in. 

The Lie algebra $M_{\geq 1}$ has a decreasing filtration
$\{M_{\geq d}\}_{d=1,\cdots, 23}$ and the adjoint action of $M_0$ respects
this filtration.
We consider one step of this filtration
\begin{equation*}
M_d=\frac{M_{\geq d}}{M_{\geq d+1}}\lhd \frac{M_{\geq 1}}{M_{\geq d+1}}
\end{equation*}
and the related Hochschild-Serre spectral sequence
\begin{equation}\label{HS-final}
E_2^{r,s}=H^r\left(\frac{M_{\geq 1}}{M_{\geq d}},H^s\left(M_d,M_{-3}\right)
\right)\Rightarrow H^{r+s}\left(\frac{M_{\geq 1}}{M_{\geq d+1}},M_{-3}\right).
\end{equation}
We fix a certain degree $d$ and we study, via the above spectral sequence, how the
truncated cohomology groups change if we pass from $d$ to $d+1$.
It is easily checked that, by homogeneity, we have 
$H^2\left(\frac{M_{\geq 1}}{M_{\geq 3}},M_{-3}\right)^{M_0}
\subset C^2(M_1,M_{-3})_{\0}=0.$
Therefore, for the rest of this section, we suppose that $d\geq 3$.
Observe also that, since $M_d$ is in the center of $M_{\geq 1}/M_{\geq d+1}$  and
$M_{-3}$ is a trivial module, then $H^s\left(M_d,M_{-3}\right)=
C^s\left(M_d,M_{-3}\right)$ and $M_{\geq 1}/M_{\geq d}$ acts trivially on it.

Consider the following exact sequence deduced from the spectral sequence 
(\ref{HS-final})
$$E_{\infty}^{1,0}=H^1\left(\frac{M_{\geq 1}}{M_{\geq d}},M_{-3}\right)\hookrightarrow
 H^1\left(\frac{M_{\geq 1}}{M_{\geq d+1}}, M_{-3}\right)\twoheadrightarrow E_{\infty}^{0,1}.$$
From the above Lemma \ref{M-commutators} (using that $d\geq 3$), 
we get that the first two terms are equal to 
$$H^1\left(\frac{M_{\geq 1}}{M_{\geq d}},M_{-3}\right)=
 H^1\left(\frac{M_{\geq 1}}{M_{\geq d+1}}, M_{-3}\right)=C^1\left(\frac{M_{\geq 1}}
{M_{\geq 3}+\langle \widetilde{x_1D_1}+\widetilde{x_2D_2}\rangle_F}, M_{-3}\right).
$$
and therefore we deduce that $E_{\infty}^{0,1}=0$. 
Together with the vanishing $E_{\infty}^{0,2}=0$ proved in the
Lemma \ref{M-inv-(0,2)} below, we deduce the following exact diagram  
$$\xymatrix{
 C^1(M_d,M_{-3}) \ar@{^{(}->}[d]_{\alpha} & &  \\
 H^2\left(\frac{M_{\geq 1}}{M_{\geq d}},M_{-3}\right) \ar@{->>}[d] & & \\
 E_{\infty}^{2,0} \ar@{^{(}->}[r] & H^2\left(\frac{M_{\geq 1}}{M_{\geq d+1}},
M_{-3}\right) \ar@{->>}[r] & E_{\infty}^{1,1}    \\
}$$
We take the cohomology with respect to $M_0$ and use the Lemmas 
\ref{M-inv-(1,1)}, \ref{M-inv-(0,1)} and \ref{M-H^1-(0,1)}.

Observe that the cocycle ${\rm inv}\circ [-,-]\in (E_{\infty}^{1,1})^{M_0}$ 
which appears for 
$d=15$ is annihilated by the differential of ${\rm inv}\in (E_2^{0,1})^{M_0}
=C^1(M_d,M_{-3})^{M_0}$ for $d=17$. Moreover the element $\langle x_iD_j\to 
\delta_{ij} {\rm inv}\rangle\in
H^1(M_0,C^1(M_{17},M_{-3}))$ does not belong to the kernel of the differential map
$$\d : H^1(M_0,C^1(M_{17},M_{-3}))=H^1(M_0,E_2^{0,1})\to H^1(M_0,E_2^{2,0}).$$
Indeed the element $\langle x_iD_j\to \delta_{ij} {\rm inv}\rangle$ does not vanish on 
$T_M$ and the same is true for its image through the map $\d$, while any coboundary 
of $H^1(M_0,E_2^{2,0})$ must vanish on $T_M$ by homogeneity.
Therefore the only cocycles that contribute to the required cohomology group are 
$\{{\rm Sq}(1), {\rm Sq}(D_1), {\rm Sq}(D_2)\}$.
\end{proof}

The remaining part of this section is devoted to prove the Lemmas that were used 
in the proof of the above Proposition. In the first Lemma, we show the 
vanishing of the term $E_{\infty}^{0,2}$ of the above spectral sequence 
(\ref{HS-final}).

\begin{lem}\label{M-inv-(0,2)}
In the above spectral sequence (\ref{HS-final}), we have $E_3^{0,2}=0$.
\end{lem}
\begin{proof}
By definition, $E_3^{0,2}$ is the kernel of the map
$$
\d: C^2\left(M_d,M_{-3}\right)=E_2^{0,2}\to E_2^{2,1}=
H^2\left(\frac{M_{\geq 1}}{M_{\geq d}},
C^1\left(M_d,M_{-3}\right)\right)
$$
that sends a $2$-cochain $f$ to the element $\d f$ given by
$\d f_{(E,F)}(G)=-f([E,F],G)$
whenever ${\rm deg}(E)+{\rm deg}(F)=d$ and $0$ otherwise.\\
The subspace of coboundaries $B^2\left(\frac{M_{\geq 1}}{M_{\geq d}},
C^1\left(M_d, M_{-3}\right)\right)$ is the image of the map
$$\partial: C^1\left(\frac{M_{\geq 1}}{M_{\geq d}}, C^1\left(M_d,
M_{-3}\right)\right)\to C^2\left(\frac{M_{\geq 1}}{M_{\geq d}}, C^1\left(M_d,
M_{-3}\right)\right)$$
that sends the element $g$ to the element $\partial g$ given by
$ \partial g_{(E,F)}(G)=-g_{[E,F]}(G).$
Hence $\partial g$ vanishes on the pairs $(E,F)$ for which ${\rm deg}(E)+{\rm deg}(F)=d$.

Therefore, if an element $f\in C^2\left(M_d,M_{-1}\right)$ is in the kernel of
$\d$, that is $\d f=\partial g$ for some $g$ as before, then it should satisfy
$f([E,F],G)=0$
for every $E, F, G$ such that ${\rm deg}(G)=d$ and ${\rm deg}(E)+{\rm deg}(F)=d$. By letting
$E$ vary in $M_1$ and $F$ in $M_{d-1}$, the bracket $[E,F]$ varies in all
$M_d$ by Lemma
\ref{M-commutators} (note that we are assuming $d\geq 3$).
Hence the preceding condition implies that $f=0$.
\end{proof}

In the newt two Lemmas, we compute the $M_0$-invariants and the first $M_0$-cohomology
group of the term $E_2^{0,1}=C^1(M_d,M_{-3})$ in the above spectral sequence (\ref{HS-final}).

\begin{lem}\label{M-inv-(0,1)}  
Define ${\rm inv}: M_{17}\to M_{-3}$ by ($i\neq j$): 
$$
{\rm inv}(x_i^3x_j^3\widetilde{D_i})=\sigma(i) D_i \: \text{ and }\:
{\rm inv}(x_i^4x_j^2\widetilde{D_i})=\sigma(i)D_j,
$$ 
where $\sigma(i)=1$ or $-1$ if $i=1$ or $i=2$, respectively. Then
$$C^1\left(M_d,M_{-3}\right)^{M_0}=\begin{cases}
\langle {\rm inv} \rangle & \text{ if } d=17,\\
0 & \text{ otherwise.} 
\end{cases}$$
\end{lem}
\begin{proof}
By homogeneity, we can assume that $d\equiv 2 \mod 5$.  We will consider the 
various cases separately.

$\fbox{ d=7 }$
Since $M_7=A(2)_3=\langle x_1^3, x_1^2x_2, x_1x_2^2, x_2^3\rangle_F $ with weights 
 $(1,3)$, $(0,4)$, $(4,0)$, $(3,1)$ respectively, a homogeneous cochain
$g\in C^1(M_d,M_{-3})^{M_0}$ can take the non-zero values 
$g(x_i^2x_j)=a_i D_j$ for $i\neq j$. We obtain the vanishing from the following 
$M_0$-invariance condition
$$0=(x_jD_i\circ g)(x_i^3)=-g(3x_i^2x_j)=-3 a_iD_j.$$

$\fbox{ d=12 }$
A homogeneous cochain $g \in C^1\left(M_{12},M_{-3}\right)_{\0}$ 
can take the non-zero values: $g(x_i^4x_jD_j)=\eta_j^i D_i$
for $i\neq j$. We get the vanishing of $g$ by mean of the
following cocycle condition
\begin{equation}\label{W-coboundary}
0=\d g_{x_iD_j}(x_i^3x_j^2D_j)=-2g(x_i^4x_jD_j)=-2\eta_j^i D_i.
\end{equation}

$\fbox{ d=17} $
$M_{17}=\widetilde{W(2)}_5=\oplus_{i\neq j, k} \langle x_i^4x_j^2 \widetilde{D_k}\rangle_F
\oplus_{i\neq j, k} \langle x_i^3x_j^3\widetilde{D_k}\rangle$ 
with weights $\epsilon_i-\epsilon_j-\epsilon_k$ and $-\epsilon_k$,
respectively. A homogeneous cochain $g\in C^1(M_d,M_{-3})^{M_0}$ 
can take the non-zero values: $g(x_i^3x_j^3\widetilde{D_i})=\alpha_i D_i$ and 
$g(x_i^4x_j^2\widetilde{D_i})=\beta_i D_j$ for $i\neq j$.
The only $M_0$-invariance conditions are the following 
$$\begin{sis}
&(x_iD_j\circ g)(x_i^2x_j^4\widetilde{D_i})=-g(4x_i^3x_j^3\widetilde{D_i}-x_i^2x_j^4
\widetilde{D_j})=
(-4\alpha_i+\beta_j)D_i,\\
&(x_iD_j\circ g)(x_i^2x_j^4\widetilde{D_j})=[x_iD_j,g(x_i^2x_j^4\widetilde{D_j})]-g(4x_i^3x_j^3
\widetilde{D_j})=(-\beta_j-4\alpha_j)D_j,\\
&(x_iD_j\circ g)(x_i^3x_j^3\widetilde{D_i})=(-\alpha_i-3\beta_i+\alpha_j)D_j,
\end{sis}$$
from which it follows that $\alpha_i=\beta_i=-\alpha_j=-\beta_j$,  that is 
$g$ is a multiple of ${\rm inv}$.

$\fbox{ d=22 }$
Since $M_{22}=A(2)_8=\langle x_1^4x_2^4\rangle$ with weight $(2,2)$,
there are no homogeneous cochains.
\end{proof}

\begin{lem}\label{M-H^1-(0,1)}
We have that
$$H^1\left(M_0,C^1\left(M_d,M_{-3}\right)\right)=
\begin{sis}
& \bigoplus_{i=1}^2 \langle \overline{{\rm Sq}(D_i)} \rangle_F & \text{ if } d=12,\\
& \langle x_iD_j\mapsto \delta_{ij} {\rm inv}\rangle_F & \text{ if } d=17,\\
& 0 &\text{ otherwise, }   
\end{sis}$$
where $\overline{{\rm Sq}(D_i)}$ is the restriction of ${\rm Sq}(D_i)$ to 
$M_0\times M_{12}$ and 
${\rm inv}:M_{17}\to M_3$ is the cochain defined in Lemma \ref{M-inv-(0,1)}. 
\end{lem}
\begin{proof}
By homogeneity, we can assume that $d\equiv 2 \mod 5$ and consider the various cases 
separately.

$\fbox{ d=7 }$
Consider a homogeneous cochain $f\in C^1(M_0,C^ 1(M_{7}, M_{-3}))_{\0}$. By applying 
cocycle conditions of the form $0=\d f_{(x_kD_k,x_iD_j)}$, it follows that $f_{x_kD_k}\in 
C^1(M_7,M_{-3})^{M_0}=0$ (by the preceding Lemma). By homogeneity, $f$ can take only  
the following non-zero values for $i\neq j$: $f_{x_iD_j}(x_j^3)=b_i D_i$ and
$f_{x_iD_j}(x_ix_j^2)=c_i D_j$.

By possibly modifying $f$ with a coboundary $\d g$ (see Lemma \ref{M-inv-(0,1)}), 
we can suppose that  $b_i=b_j=0$.  The following cocycle condition
$$0=\d f_{(x_iD_j,x_jD_i)}(x_i^2x_j)=[x_iD_j,f_{x_jD_i}(x_jx_i^2)]-f_{x_jD_i}(x_i^3)+
f_{x_iD_j}(2x_ix_j^2)=$$
$$=(-c_j-b_j+2c_i)D_j=(-c_j+2c_i)D_j,$$
gives that $c_j=2c_i=4c_j$ and hence that $c_i=c_j=0$, that is $f=0$.

$\fbox{ d=12} $
First of all, the cocycles  $\overline{\Sq(D_i)}$ for $i=1, 2$ belong
to $H^1(M_0,C^1(M_{12},$ $  M_{-3}))_{\0}$ since they are restriction 
of global cocycles and, moreover, they are independent
as it follows from the formulas (\ref{indep3}) and (\ref{indep3bis}). 
It remains to show that the above cohomology group has dimension 
less than or equal to $2$.
Consider a homogeneous cocycle $f\in C^1(M_0, C^1(M_{12}, M_{-3}))_{\0}$.  
First of all, we observe that $f$ must satisfy $f_{x_iD_i}=0.$
Indeed, by the formula (\ref{restriction}), we have $0=\d f_{| x_iD_i}=x_iD_i\circ f -
\d(f_{| x_iD_i})$ from which, since the first term vanishes for homogeneity reasons,
it follows that
$f_{| x_iD_i}\in C^1\left(M_{p-1},M_{-1}\right)^{M_0}$ which is zero by
the previous Lemma \ref{M-inv-(0,1)}.  Therefore $f$ can take only the following
non-zero values (for $i\neq  j$):
$$\begin{sis}
&f_{x_iD_j}(x_i^{3}x_j^2D_j)=\alpha_{ij} D_i, \\
&f_{x_iD_j}(x_i^{4}x_jD_i)=\beta_{ij} D_i,  \\
&f_{x_iD_j}(x_i^{4}x_jD_j)=\gamma_{ij}D_j.
\end{sis}$$
By possibly modifying $f$ with a coboundary (see formula (\ref{W-coboundary})),
we can assume that $\underline{\alpha_{ij}=0}$. 
The coefficients $\underline{\beta_{ij}}$ are
determined by the coefficients $\gamma_{ij}$ in virtue of 
the following cocycle condition:
$$ 0=\d f_{(x_iD_j,x_jD_i)}(x_i^{4}x_jD_j)=f_{x_iD_j}(-x_i^{3}x_j^2D_j)+f_{x_iD_j}(-x_i^{4}x_jD_i)+$$
\begin{equation*}
-[x_jD_i,f_{x_iD_j}(x_i^{4}x_jD_j)]=[-\alpha_{ij}-\beta_{ij}+
\gamma_{ij}] D_i=[-\beta_{ij}+ \gamma_{ij}] D_i.  
\end{equation*}
Therefore the cohomology group depends on the two parameters $\gamma_{12}$ and 
$\gamma_{21}$ and hence has dimension less than or equal to $2$. 


$\fbox{ d=17 }$
Consider a homogeneous cochain $f \in C^1(M_0,C^1(M_{17},M_{-3}))_{\0}$.
By imposing cocycle conditions of the form $0=\d f_{(x_kD_k,x_iD_j)}$, one obtains
that $f_{x_kD_k}\in C^1(M_{17}, M_{-3})^{M_0}=\langle {\rm inv}\rangle$ (by the 
preceding Lemma). Put $f_{x_kD_k}=\mu_k \cdot {\rm inv}$. By homogeneity, the other non-zero 
values of $f$ are (for $i\neq j$):
$$f_{x_iD_j}(x_i^2x_j^4\widetilde{D_i})=\nu_i D_i, \hspace{0,3cm} 
f_{x_iD_j}(x_i^2x_j^4\widetilde{D_j})=\tau_i D_j,\hspace{0,3cm}
f_{x_iD_j}(x_i^3x_j^3 \widetilde{D_i})=\sigma_i D_j.$$
By possibly modifying $f$ with a coboundary $\d g$, we can assume that 
$f_{x_1D_2}=0$ (see Lemma \ref{M-inv-(0,1)}).
Considering all the cocycle conditions of the form $0=\d f_{(x_1D_2,x_2D_1)}=(x_1D_2\circ 
f_{x_2D_1})-f_{x_1D_1}+f_{x_2D_2}$, one gets
$$\begin{sis}
&0=df_{(x_1D_2,x_2D_1)}(x_1^3x_2^3\widetilde{D_1})=-3\tau_2+\sigma_2-(\mu_1-\mu_2),\\
&0=df_{(x_1D_2,x_2D_1)}(x_1^3x_2^3\widetilde{D_2})=\sigma_2+3\nu_2-(\mu_1-\mu_2), \\
&0=df_{(x_1D_2,x_2D_1)}(x_1^2x_2^4\widetilde{D_2})=\sigma_2+(\mu_1-\mu_2),   \\
&0=df_{(x_1D_2,x_2D_1)}(x_1^4x_2^2\widetilde{D_1})=\tau_2-\nu_2+(\mu_1-\mu_2).
\end{sis}$$
Since the $4\times 4$ matrix 
associated to the preceding system of $4$ equations in the $4$ variables $\nu_2$, $\tau_2$,
$\sigma_2$, $\mu_1-\mu_2$ is invertible (it has determinant equal to $2$),  we conclude 
that $f_{x_2D_1}=0$ and $f_{x_1D_1}=f_{x_2D_2}$ is a multiple of ${\rm inv}$.

$\fbox{ d=22 }$
There are no homogeneous cochains.

\end{proof}

In the last Lemma, we compute the $M_0$-invariants of the term $E_{\infty}^{1,1}$
of the above spectral sequence (\ref{HS-final}).
In view of Lemma \ref{M-commutators} and the 
hypothesis $d\geq 3$, we have that 
$$E_2^{1,1}=C^1\left(\frac{M_1\oplus M_2}{\langle x_1\widetilde{D_1}+x_2\widetilde{D_2}\rangle }
\times M_d,M_{-3}\right).$$

\begin{lem}\label{M-inv-(1,1)}
In the above spectral sequence (\ref{HS-final}), we have that
$$(E_{\infty}^{1,1})^{M_0}=
\begin{cases}
\langle \overline{{\rm Sq}(1)} \rangle_F & \text{ if } d=6,\\
\langle {\rm inv}\circ [-,-] \rangle_F & \text{ if } d=15,\\
0 & \text{ otherwise, }
\end{cases}$$
where $\overline{{\rm Sq}(1)}$ is the restriction of ${\rm Sq}(1)$ to $M_1\times M_6$
and ${\rm inv}\circ [-,-]$ is defined by 
$$({\rm inv}\circ [-,-])(E,F)={\rm inv}([E,F]),$$
where ${\rm inv}:M_{17}\to M_{-3}$ is the $M_0$-invariant map defined in Lemma
\ref{M-inv-(0,1)}. 
\end{lem}
\begin{proof}
The term $(E_{\infty}^{1,1})^{M_0}=(E_3^{1,1})^{M_0}$ is the kernel of the map
$$\d:(E_2^{1,1})^{M_0}\longrightarrow (E_2^{3,0})^{M_0}\hookrightarrow 
E_2^{3,0}=H^3\left(\frac{M_{\geq 1}}{M_{\geq d}},M_{-3}\right).$$
In order to avoid confusion, during this proof we denote
with $\partial f\in B^3\left(\frac{M_{\geq 1}}{M_{\geq d}},M_{-3}\right)$ 
(instead of the usual $\d f$)
the coboundary of an element $f\in C^2\left(\frac{M_{\geq 1}}{M_{\geq d}},M_{-3}\right)$.

Note that $M_1=\langle x_1, x_2\rangle_F$ with weights respectively  $(4,3)$ and $(3, 4)$ 
while $M_2/ $ $( x_1\widetilde{D_1}+x_2\widetilde{D_2})=\langle 
x_1\widetilde{D_1}=- x_2\widetilde{D_2}, x_1\widetilde{D_2}, x_2\widetilde{D_1}\rangle_F $ 
with weights respectively $(2,2)$, $(3,1)$ and $(1,3)$. Note also that after the identification 
$x_1\widetilde{D_1}=-x_2\widetilde{D_2}$, the action of 
$M_0=W(2)_0$ on $M_2/(x_1\widetilde{D_1}+x_2\widetilde{D_2})$ is given by (for $i\neq j$) 
$[x_iD_j,x_i\widetilde{D_j}]=0$, $[x_i D_j,x_1\widetilde{D_1}]=\sigma(j)x_i\widetilde{D_j}$ 
and $[x_iD_j, x_j\widetilde{D_i}]=-2\sigma(j)x_1\widetilde{D_1}$ where $\sigma(j)=1,-1$ 
if respectively $j=1, 2$.

By homogeneity, an $M_0$-invariant cochain of $E_2^{1,1}$ 
can assume non-zero values only on $M_1\times M_d$ if 
$d\equiv 1 \mod 5$ or on $M_2\times M_d$ if $d\equiv 0 \mod 5$.  We will consider the various 
cases separately.

$\fbox{ d=5 }$
$M_5=\widetilde{W(2)}_1=\oplus_{i,k}\langle x_i^2 \widetilde{D_k}\rangle_F
\oplus_{k} \langle x_1x_2\widetilde{D_k}\rangle_F$ 
with weights $(2,2)+2\epsilon_i-\epsilon_k$ and $(3,3)-\epsilon_k$, respectively.
A homogeneous cochain $g\in (E_2^{1,1})^{M_0}$ 
can take only the following non-zero values (for $i\neq j$):
$$\begin{sis}
 &g(x_1\widetilde{D_1}, x_ix_j \widetilde{D_i})=h_i D_i, \hspace{4,8cm}
g(x_1\widetilde{D_1}, x_j^2 \widetilde{D_j})=k_i D_i,\\
 & g(x_i\widetilde{D_j}, x_j^2\widetilde{D_i})=l_i D_i, \hspace{0,7cm} 
 g(x_i\widetilde{D_j}, x_ix_j\widetilde{D_i})=m_i D_j, \hspace{0,7cm}
 g(x_i\widetilde{D_j}, x_j^2\widetilde{D_j})=n_i D_j . 
\end{sis}$$
Consider the following $M_0$-invariance conditions:
$$\begin{sis}
&0=(x_iD_j\circ g)(x_j\widetilde{D_i}, x_j^2 \widetilde{D_j})=2(\sigma(j)k_i-m_j)D_i,\\ 
&0=(x_iD_j\circ g)( x_j\widetilde{D_i},x_ix_j\widetilde{D_i})=(2\sigma(j)h_i-n_j+m_j)D_i,\\
&0=(x_iD_j\circ g)(x_i\widetilde{D_j}, x_j^2\widetilde{D_i})=(-l_i-2m_i+n_i)D_j,\\
&0=(x_iD_j\circ g)(x_j\widetilde{D_i}, x_ix_j \widetilde{D_j})=(-m_j+2\sigma(j)h_j-l_j)D_j,\\
&0=(x_iD_j\circ g)(x_j\widetilde{D_i},x_i^2\widetilde{D_i})=(-n_j+2\sigma(j)k_j+l_j)D_j.
\end{sis}$$
From the first $3$ equations one obtains that $(m_j, n_j, l_j)=\sigma(j)(k_i, 2h_i+k_i, 
2h_i-k_i)$
and substituting in the last two equations one finds that $h_i=h_j:=h$ and $k_i=k_j:=k$. 
 
Suppose now that $g$ is in the kernel of the map $\d$, that is $\d g=-\partial f$ for some
$f\in C^2(M_{\geq 1}/M_{\geq 5},M_{-3})$. Applying $0=\d g+\partial f$ to the triples 
$(x_i,x_j\widetilde{D_j},x_ix_j)$ and $(x_j,x_i\widetilde{D_j},x_ix_j)$ for $i\neq j$,
we get the two conditions
$$\begin{sis}
& f(x_ix_j\widetilde{D_j},x_ix_j)=2 \sigma(i) g(x_j \widetilde{D_j},x_i^2 \widetilde{D_i})=
-2kD_j,\\
& f(x_ix_j\widetilde{D_j},x_ix_j)=2\sigma(j) g(x_i\widetilde{D_j}, x_j^2\widetilde{D_j})=
-2(2h+k)D_j,
\end{sis}$$
from which it follows that $h=0$. Considering now the triples $(x_i,x_r \widetilde{D_r}, 
x_ix_j)$ and $(x_i,x_r\widetilde{D_r}, x_j^2)$ (for $i\neq j$ and some $r=1,2$), we get
$$\begin{sis}
& f(x_ix_r\widetilde{D_r},x_ix_j)=2 \sigma(i) g(x_r \widetilde{D_r}, x_i^2\widetilde{D_i})=
2\sigma(i)\sigma(r) k D_j,\\
& f(x_ix_r\widetilde{D_r}, x_j^2)=-\sigma(i)g(x_r\widetilde{D_r}, x_ix_j\widetilde{D_i}-
2 x_j^2 \widetilde{D_j})= 2 \sigma(i)\sigma(r) k D_i. 
\end{sis}$$
Finally, considering the triple $(x_i,x_i^2 D_i, x_ix_j D_i)$ and using the two preceding
relations, we get
$$0=(\d g+\partial f)(x_i,x_i^2 D_i, x_ix_j D_i)=
f([x_i^2D_i, x_i],x_ix_jD_i)-f([x_ix_jD_i,x_i],x_i^2D_i)=$$
$$=3 f(x_jx_i D_i, x_i^2)-f(x_i^2 D_i, x_i x_j)=-6kD_j-2k D_j=2 kD_j,$$
from which it follows that $k=0$, that is $g=0$.   

$\fbox{ d=6 }$
First of all, the cocycle $\overline{{\rm Sq}(1)}$ is an element of 
$(E_{\infty}^{1,1})^{M_0}$ since it is the restriction of a global cocycle and
is non-zero in virtue of formulas (\ref{indep2}) and (\ref{indep2bis}).
Therefore, it remains to show that the dimension of $(E_{\infty}^{1,1})^{M_0}$ 
is at most $1$. Consider a cochain $g\in (E_2^{1,1})^{M_0}$, which is in particular 
homogeneous. Since $M_6=W(2)_2=\oplus_{i,k}\langle x_i^2 D_k\rangle_F \oplus_{k} \langle x_1x_2 D_k \rangle_F$ with weights $2\epsilon_i-\epsilon_k$
and $(1,1)-\epsilon_k$ respectively, $g$ can take only 
the following non-zero values 
(for $i\neq j$):
$$\begin{sis}
&g(x_i, x_ix_j^2D_i)= c_i D_i, \hspace{1cm} &  g(x_i, x_i^2 x_j D_i)= e_i D_j,\\
&g(x_i, x_j^3 D_j)= d_i D_i, & g(x_i, x_ix_j^2 D_j)= f_i D_j.\\
\end{sis}$$ 
Consider the following $M_0$-invariance conditions:
$$\begin{sis}
& 0=(x_iD_j\circ g)(x_i, x_j^3 D_i)=-g(x_i, 3x_ix_j^2 D_i-x_j^3 D_j)=(-3c_i+d_i) D_i,\\
& 0=(x_iD_j\circ g)(x_i, x_i x_j^2 D_i)=(-c_i-2e_i+f_i)D_j,\\
& 0=(x_iD_j\circ g)(x_i, x_j^3 D_j)=[x_iD_j,d_iD_i]-g(x_i,3x_ix_j^2 D_j)=(-d_i-3f_i)D_i,\\
& 0=(x_iD_j\circ g)(x_j, x_i x_j^2 D_i)=(-c_i-2f_j+e_j)D_i.\\
\end{sis}$$
From the first $3$ equations, one gets $(f_i, c_i, d_i)=(e_i, -e_i, 2e_i)$ and substituting 
into the last equation  one finds $e_i=e_j:=e$. Therefore
$(E_{\infty}^{1,1})^{M_0}$ depends on one parameter and hence has 
dimension at most $1$.
  
$\fbox{ d=10 }$
$M_{10}=A(2)_4=\oplus_{i}\langle x_i^4\rangle_F \oplus_{i\neq j} \langle
x_i^3x_j\rangle_F \oplus_{i\neq j} \langle x_i^2x_j^2\rangle_F $ with weights 
$2\epsilon_i-2\epsilon_j$, $\epsilon_i-\epsilon_j$ and $(0,0)$, respectively.
Consider a cochain $g\in (E_2^{1,1})^{M_0}$. 
 Since $M_0$ acts transitively on $M_2$, to prove the vanishing of $g$
it is enough to prove that $g(x_1\widetilde{D_2}, -)=0.$ By homogeneity the 
only possible non-zero 
such values are $g(x_1\widetilde{D_2}, x_1^3x_2)$ and $g(x_1\widetilde{D_2}, x_1^4)$ and 
the vanishing follows from the $M_0$-invariance conditions:
$$\begin{sis}
&0=(x_1D_2\circ g)(x_1\widetilde{D_2}, x_1^2x_2^2)=-g(x_1\widetilde{D_2},2 x_1^3x_2),\\
&0=(x_1\widetilde{D_2}\circ g)( x_1\widetilde{D_2},x_1^3 x_2)=
[x_1D_2, g(x_1\widetilde{D_2}, x_1^3x_2)]-g(x_1\widetilde{D_2},x_1^4). 
\end{sis}$$

$\fbox{ d=15 }$
$M_{15}=W(2)_5=\oplus_{i\neq j, k}\langle x_i^4 x_j^2 D_k\rangle_F \oplus_{i\neq j, k}
\langle x_i^3 x_j^3 D_k\rangle_F$ with weights 
$-\epsilon_i+2\epsilon_j-\epsilon_k$ and $3\epsilon_i+3\epsilon_j-\epsilon_k$, respectively.
A homogeneous cochain $g\in (E_2^{1,1})^{M_0}$
can take only the following non-zero values (for $i\neq j$):
$$\begin{sis}
&g(x_1\widetilde{D_1}, x_i^3x_j^3 D_i)=p_i D_i, & 
g(x_1\widetilde{D_1}, x_i^2x_j^4 D_j)=q_i D_i,\\
 & g(x_i\widetilde{D_j}, x_i^2x_j^4 D_i)=r_i D_i, \hspace{0.3cm}
 g(x_i\widetilde{D_j}, x_i^3x_j^3 D_i)=s_i D_j, & 
 g(x_i\widetilde{D_j}, x_i^2 x_j^4 D_j)=t_i D_j . 
\end{sis}$$
Consider the following $M_0$-invariance conditions:
$$\begin{sis}
& 0=(x_iD_j\circ g)(x_j\widetilde{D_i}, x_i^2x_j^4D_j)=2(\sigma(j)q_i-2s_j)D_i,\\
& 0=(x_iD_j\circ g)(x_j\widetilde{D_i}, x_i^3x_j^3D_i)=(2\sigma(j)p_i-3t_j+s_j)D_i,\\
& 0=(x_iD_j\circ g)(x_i\widetilde{D_j}, x_i^2x_j^4D_i)=(-r_i+s_i+t_i)D_j,\\
& 0=(x_iD_j\circ g)(x_j\widetilde{D_i}, x_i^3x_j^3D_j)=(-s_j+2\sigma(j)p_j-3r_j)D_j,\\
& 0=(x_iD_j\circ g)(x_j\widetilde{D_i}, x_i^4x_j^2D_i)=(-t_j+2\sigma(j)q_j+r_j)D_j.
\end{sis}$$
From the first $3$ equations one gets that $(s_j,t_j, r_j)=\sigma(j)(-2q_i,-p_i+q_i, -p_i-q_i)$ 
and substituting in the last two equations one obtains $p_i=p_j:=p$ and $q_i=q_j:=q$.

Suppose now that $g$ is in the kernel of $\d$, that is $\d g=-\partial f$ for some 
$f\in C^2(M_{\geq 1}/M_{\geq 15},M_{-3})$. Applying $0=\d g+\partial f$
to the triple $(x_i, x_j \widetilde{D_j}, x_i^4x_j \widetilde{D_i})$ for $i\neq j$, we get
$$0=-f([x_i,x_j\widetilde{D_i}],x_i^4 x_j \widetilde{D_i})+g([x_i,x_i^4x_j \widetilde{D_i}],
x_j\widetilde{D_i})=f(x_ix_jD_i, x_i^4x_j \widetilde{D_i}).$$
Considering the triple $(x_j,x_i\widetilde{D_i},x_i^4x_j 
\widetilde{D_i})$ and using the preceding vanishing,  
 we obtain $$0=-f(x_ix_j\widetilde{D_i},x_i^4x_j \widetilde{D_i})
-g(x_i\widetilde{D_i},x_i^4x_j^2\widetilde{D_i})=-\sigma(i)q D_j,$$
that is $q=0$. For $p=1$, one obtains that $g={\rm inv}\circ [-, -]$, which is clearly
in the kernel of $\d$ since it is the restriction of a cocycle on $M_{\geq 1}\times M_{\geq 1}$.

$\fbox{ d=16 }$
$M_{16}=A(2)_6=\oplus_{i\neq j}\langle x_i^2 x_j^4\rangle_F\oplus 
\langle x_1^3 x_2^3\rangle_F$ 
with weights $2\epsilon_j$ and $(1,1)$, respectively.
A homogeneous cochain $g\in (E_2^{1,1})^{M_0}$ 
can take the non-zero values (for $i\neq j$):
$g(x_i, x_i^2x_j^4)=a_i D_i$ and  $g(x_i, x_1^3x_1^3)=b_iD_j $ for $i\neq j$.
We get that $a_i=a_j=b_i=b_j:=a$ by considering the following $M_0$-invariant conditions:
$$\begin{sis}
&0=(x_iD_j\circ g)(x_i, x_i^2x_j^4)=[x_iD_j,a_iD_i]-g(x_i,4x_i^3x_j^3)=(-a_i+b_i)D_j,\\
&0=(x_iD_j\circ g)(x_j, x_i^2x_j^4)=-g(x_i,x_i^2x_j^4)-g(x_j,4x_i^3x_j^3)=(-a_i+b_j)D_i.\\
\end{sis}$$
Now suppose that $g$ is in the kernel of the map $\d$, that is $\d g=-\partial f$ for some
$f\in C^2(M_{\geq 1}/M_{\geq 16},M_{-3})$.
Applying the relation $0=\d g+\partial f$ to the two triples $(x_i,x_j,x_i^4x_j^2 D_i)$
and $(x_i,x_j,x_i^3x_j^3D_i)$ for $i\neq j$, we get
$$\begin{sis}
& f(x_1\widetilde{D_1}+x_2\widetilde{D_2},x_i^4x_j^2D_i)=g([x_i,x_i^4x_j^2D_i],x_j)-
g([x_j,x_i^4x_j^2D_i],x_i)=a,\\
& f(x_1\widetilde{D_1}+x_2\widetilde{D_2},x_i^3x_j^3D_i)= g([x_i,x_i^3x_j^3D_i],x_j)-
g([x_j,x_i^3x_j^3D_i],x_i)=a. 
\end{sis} \hspace{0,7cm} (*)$$  
Considering the triples $(x_1,x_1\widetilde{D_1}+x_2\widetilde{D_2}, 
x_1^3x_2^2\widetilde{D_1})$ 
and $(x_1,x_1\widetilde{D_1}+x_2\widetilde{D_2}, x_1^2x_2^3\widetilde{D_2})$ and 
using $(*)$, we get 
$$\begin{sis}
& f(x_1^2D_1+x_1x_2D_2,x_1^3x_2^2\widetilde{D_1})=-f(x_1\widetilde{D_1}+x_2\widetilde{D_2},
x_1^4x_2^2D_1)-g(x_1,x_1^3x_2^3)=-2a,\\
& f(x_1^2D_1+x_1x_2D_2,x_1^2x_2^3\widetilde{D_2})=-f(x_1\widetilde{D_1}+x_2\widetilde{D_2},
x_1^3x_2^3D_2)+g(x_1,x_1^3x_2^3)=0. 
\end{sis} (**)$$
Finally, using the relations $(**)$, we get the vanishing of $g$ by mean of the following
$$0=(\d g+\partial f)(x_1,x_1^2D_1+x_1x_2D_2,x_1^2x_2^3)=-f([x_1,x_1^2D_1+x_1x_2D_2],
x_1^2x_2^3)+$$
$$+f([x_1,x_1^2x_2^3],x_1^2D_1+x_1x_2D_2)-g([x_1^2D_1+x_1x_2D_2,x_1^2x_2^3],x_1)=$$
$$=f(-x_1^3x_2^2\widetilde{D_1}+2x_1^2x_2^3\widetilde{D_2},x_1^2D_1+x_1x_2D_2)-
g(-x_1^3x_2^3, x_1)=-2a-a=-3 a.$$

$\fbox{ d=20 }$
$M_{20}=\widetilde{W(2)}_6=\oplus_{i\neq j, k} \langle x_i^4x_j^3\widetilde{D_k}\rangle_F$
with weights $\epsilon_i-\epsilon_k$. A homogeneous cochain $g\in (E_2^{1,1})^{M_0}$ 
can take only the non-zero 
values $g(x_i\widetilde{D_j}, x_i^4x_j^3\widetilde{D_j})=\lambda_i D_i$. The vanishing follows
from the $M_0$-invariance conditions
$$0=(x_iD_j\circ g)(x_i\widetilde{D_j}, x_i^3 x_j^4\widetilde{D_j})=-4 g(x_i\widetilde{D_j}, 
x_i^4x_j^3\widetilde{D_j})=\lambda_i D_i.$$

$\fbox{ d=21 }$
Since $M_{21}=W(2)_7=\langle x_1^4x_2^4D_1, x_1^4x_2^4D_2\rangle_F$ with weights 
$(3,4)$ and $(4,3)$ respectively, there are no homogeneous cochains.

\end{proof}


\begin{thebibliography}{TESI}

\addcontentsline{toc}{chapter}{References}


\bibitem[BW88]{BW} Block R.E., Wilson R.L. Classification of the restricted simple
Lie algebras.  \emph{J. Algebra}  114  (1988): 115--259.

\bibitem[CE48]{CE} Chevalley C., Eilenberg S. Cohomology theory of Lie groups
and Lie algebras. \emph{Trans. Amer. Math. Soc.} 63 (1948): 85--124.

\bibitem[GER64]{GER1} Gerstenhaber M. On the deformation of rings and algebras.
\emph{Ann. of Math.} 79 (1964): 59-103. 

\bibitem[HS53]{HS} Hochschild G., Serre J.P. Cohomology of Lie algebras.
\emph{Ann. of Math.} 57 (1953): 591--603.

\bibitem[KUT90]{KUT1} Kuznetsov. M. I.
On Lie algebras of contact type. \emph{Commun. Algebra} 18 (1990): 2943-3013.

 
\bibitem[MEL80]{MEL} Melikian G.M. Simple Lie algebras of characteristic $5$
(Russian).  \emph{Uspekhi Mat. Nauk} 35 (1980): 203--204.

\bibitem[PS01]{PS3} Premet A., Strade H. Simple Lie algebras of small characteristic.
III: The toral rank 2 case. \emph{J. Algebra}  242 (2001): 236--337.

\bibitem[RUD71]{RUD} Rudakov A. N. Deformations of simple Lie algebras (Russian).
 \emph{Izv. Akad. Nauk SSSR Ser. Mat.} 35 (1971): 1113--1119.

\bibitem[STR04]{STR} Strade H. Simple Lie algebras over fields of positive 
characteristic I: Structure theory. \emph{De Gruyter Expositions in Mathematics} 38 
(2004), Walter de Gruyter, Berlin.

\bibitem[VIV1]{VIV1} Viviani F. Deformations of restricted simple Lie algebras I.
math.RA/0612861.

\bibitem[VIV2]{VIV2} Viviani F. Deformations of restricted simple Lie algebras II.
math.RA/0702499.





\end{thebibliography}
\end{document}